\documentclass[dvipsnames,a4paper,twoside]{article}
\usepackage{amsthm,amsmath,amssymb,xcolor,soul,xspace,graphicx}
\usepackage[bookmarks=true,bookmarksopen=true]{hyperref}

\usepackage[ruled,algo2e]{algorithm2e}

\theoremstyle{plain}
\newtheorem{theorem}{Theorem}[section]

\newtheorem{corollary}[theorem]{Corollary}

\theoremstyle{definition}
\newtheorem{definition}[theorem]{Definition}

\newtheorem{assumption}[theorem]{Assumption}

\theoremstyle{remark}
\newtheorem{remark}[theorem]{Remark}

 \numberwithin{equation}{section}

\title{Bilinear control of semilinear elliptic PDEs: Convergence of a semismooth Newton method\thanks{The first and third authors were supported by MCIN/ AEI/10.13039/501100011033/ under research project PID2020-114837GB-I00. The second author was supported by the Hellenic Foundation for Research and Innovation (H.F.R.I.) under the ``First Call for H.F.R.I. Research Projects to support Faculty members and Researchers and the procurement of high-cost research equipment grant" (Project Number: 3270)}}

\author{Eduardo Casas\thanks{Departamento de Matem\'{a}tica Aplicada y Ciencias de la Computaci\'{o}n, E.T.S.I. Industriales y de Telecomunicaci\'on, Universidad de Cantabria, 39005 Santander, Spain
(\texttt{eduardo.casas@unican.es})},
\and
Konstantinos Chrysafinos\thanks{Department of Mathematics, School of Applied Mathematics and Physical Sciences, National Technical University of Athens, Zografou Campus, Athens 15780, Greece and IACM, FORTH, 20013 Heraklion, Crete, Greece. (\texttt{chrysafinos@math.ntua.gr})}
\and
Mariano Mateos\thanks{Departamento de Matem\'{a}ticas, Campus de Gij\'on, Universidad de Oviedo, 33203, Gij\'on, Spain (\texttt{mmateos@uniovi.es})}
}

\ifpdf
\hypersetup{
  pdftitle={Bilinear control of semilinear elliptic PDEs: Convergence of a semismooth Newton method},
  pdfauthor={E.~Casas,  K. Chrysafinos, and M.~Mateos}
}
\fi

\newcommand{\dx}{\,\mathrm{d}x}

\newcommand{\Pb}{\mbox{\rm (P)}\xspace}

\newcommand{\uad}{U_{\rm ad}}

\newcommand{\proj}{\operatorname{Proj}}

\newcommand{\A}{\mathcal{A}}
\newcommand{\ci}{\chi_{{}_{\mathbb{I}_u}}\!}
\newcommand{\ca}{\chi_{{}_{\mathbb{A}_u}}\!}

\newcommand{\dimension}{d}
\newcommand{\conormal}{n\!}
\newcommand{\tichonov}{\nu}

\newcommand{\umin}{{\alpha}}
\newcommand{\umax}{{\beta}}

\begin{document}

\maketitle

\begin{abstract}
In this paper, we carry out the analysis of the semismooth Newton method for bilinear control problems related to semilinear elliptic PDEs.  We prove existence, uniqueness and regularity for the solution of the state equation, as well as differentiability properties of the control to state mapping. Then, first and second order optimality conditions are obtained. Finally, we prove the superlinear convergence of the semismooth
Newton method to local solutions satisfying no-gap second order sufficient optimality conditions as well as a strict complementarity condition.
\end{abstract}

\begin{quote}
\textbf{Keywords:}
optimal control,  bilinear control, semilinear elliptic equations, semismooth Newton method
\end{quote}

\begin{quote}
\textbf{AMS Subject classification: }
35J61,  
49K20,  
49M15,  
49M05 
\end{quote}

\section{Introduction} \label{S1}
In this paper, we propose a semismooth Newton method to solve the following bilinear optimal control problem:
\[ \Pb \min_{u \in \uad} J(u) :=  \int_\Omega L(x,y_u(x)) \dx + \frac{\tichonov}{2} \int_\Omega u^2(x) \dx, \]
where $y_u$ is the state associated with the control $u$ solution of

\begin{equation}
\left\{\begin{array}{l} Ay + a(x,y) + uy= 0\ \  \mbox{in } \Omega,\vspace{2mm}\\  \partial_{\conormal_A} y = g\ \ \mbox{on }\Gamma. \end{array}\right.
\label{E1.1}
\end{equation}
Here $\Omega \subset \mathbb{R}^{\dimension}$, $\dimension =2$ or $3$, is a bounded open connected set with a Lipschitz boundary $\Gamma.$ The precise assumptions of the problem will be given in Sections \ref{S2} and \ref{S3}. For the moment, we underline that the parameter $\tichonov >0$ and the admissible set of controls is defined as
\[ \uad = \{ u \in L^2(\Omega) :  \umin \leq u(x) \leq \umax \text{ a.e. in }  \Omega \}. \]

Our main goal is to prove convergence of a semi-smooth Newton method for the bilinear control problem $\Pb.$ Bilinear controls have numerous applications in biology, ecology, socio-economy and engineering; see \cite{Bruni-DiPillo-Koch1974,Mohler1973,GlowinskiSongYuanYue2022}. The key structural feature of such problems is the "bilinear" structure of the control; that is the nonlinear multiplicative coupling $uy$ of the control variable $u$ to its state variable $y$, in contrast to the classical optimal control setting (see for instance \cite{Troltzsch2010}) where the control typically appears in a linear way at the right hand side.

Semi-smooth Newton type methods are well known for their computational effectivity and their robust performance in a variety of optimization problems; see \cite{Hintermuller-Ito-Kunisch03,Ulbrich11} and references within. For various results related to the use of semi-smooth type methods within the context of PDE-constrained optimization we refer the reader to the books of \cite{Hinze-Pinnau-Ulbrich-Ulbrich2009,Ito-Kunisch2008,Ulbrich11} and references within.

Despite its wide applicability, results regarding convergence properties of semi-smooth methods associated to nonlinear PDE constrained optimization problems are very limited; see \cite{Amstutz-Laurain,Casas-Mateos2023,Mannel-Rund2021,Pieper2015}. In \cite{Amstutz-Laurain} a condition implying the convexity of the associated functional is used, while in \cite{Mannel-Rund2021,Pieper2015} convergence is shown under a strict second order condition assumption, of the form $J''({\bar u})v^2 \geq \kappa \|v\|^2_{L^2(\Omega)}$, for all $v \in L^2(\Omega)$, with $\kappa >0,$ and $\bar u$ denoting a local optimal solution, which also implies local convexity of $J$ around $\bar u$. In the recent work of \cite{Casas-Mateos2023} convergence is proved without imposing additional assumptions to the classical ones used in the finite dimensional case. In particular, under the standard assumptions of no-gap second order optimality conditions and the typical strict complementarity condition, convergence of a semi-smooth Newton method is proved under a monotonicity assumption on the semilinear term.

According to our best knowlegde the important case of bilinear controls has not being considered before. The case of bilinear controls posses additional challenges. For instance, the control enters to the PDE in a multiplicative way, and the sign of the bilinear term $uy$ is not necessarily strictly positive. As a consequence, the derivation of a suitable second order conditions substantially differs from the classical case, since various results regarding the well-posedness and differentiability properties of the control to state and adjoint-state mappings are non standard.

Our paper fills this gap in the case of bilinear controls for optimal control problems related to semi-linear elliptic PDEs. In particular, we prove the superlinear convergence of the semi-smooth Newton method, under the standard assumptions of the no-gap second order optimality conditions and a strict complementarity conditions (similar to the assumptions of \cite{Casas-Mateos2023} and to the finite dimensional case \cite{Nocedal-Wright1999}). The key ingredient of the proof is the development of a suitable second order condition on an extended cone. For the later, we prove various local well posedness and differentiability results for the associated control to state and adjoint state mappings. The second order condition allows to prove the uniform boundedness of certain generalized derivatives of the solution operator equation associated to the semi-smooth Newton method, which is a necessary result in order to exploit the abstract convergence framework of \cite{Ulbrich11}.

The remaining of the paper is organized as follows. In Section 2, we present the analysis of the state equation and, in particular, well-posedness results related to the control to state mapping. In Section 3, we study the optimal control problem, and in particular we prove  first and second order conditions for a local minimizer. In Section 4, we employ the functional framework of \cite{Ulbrich11} to study the convergence of the semi-smooth Newton method while in Section 5 we present a numerical example that verifies our theoretical findings.

\section{Analysis of the state equation} \label{S2}
\setcounter{equation}{0}
In this section we prove existence and uniqueness of \eqref{E1.1} as well as differentiablity properties of the relation control to state.
To this end, we make the following assumptions.
\begin{assumption} \label{A2.1}
The operator $A$ is defined in $\Omega$ by the expression
\[ Ay = - \sum_{i,j=1}^{\dimension} \partial_{x_j} [ a_{ij}(x) \partial_{x_i} y] \]
with $a_{ij} \in L^{\infty}(\Omega)$ for $1\leq i,j \leq \dimension$ satisfying for some $M_A,\Lambda_A >0$
\[ M_A | \xi |^2 \geq \sum_{i,j=1}^{\dimension} a_{ij}(x) \xi_i \xi_j \geq \Lambda_A | \xi |^2 \quad \text{for a.a.} \ x \in \Omega \ \text{and} \ \forall \xi \in \mathbb R^n.
\]
\end{assumption}

\begin{assumption} \label{A2.2}
We assume that $a:\Omega \times \mathbb{R} \longrightarrow \mathbb{R}$ is a Carath\'eodory function of class $C^2$ with respect to the second variable satisfying the following properties for a.a. $x \in \Omega$:
\begin{align*}
&\bullet a(\cdot,0) \in L^p(\Omega) \text{ for some } p  > \frac{d}{2},  \\
&\bullet  \frac{\partial a}{\partial y}(x,y) \geq a_0 \forall y\in\mathbb{R} \text{ for some } a_0 \geq 0, \\
&\bullet \forall M > 0 \ \exists C_{a,M} \text{ such that } \sum_{j=1}^2 \Big| \frac{\partial^j a}{\partial y^j} (x,y) \Big| \leq C_{a,M}  \  \forall | y | \leq M, \\
&\bullet \forall \varepsilon >0 \text{ and } \forall M>0 \ \exists \rho >0 \text{ such that } \ \Big|    \frac{\partial^2 a}{\partial y^2}  (x,y_1) -  \frac{\partial^2 a}{\partial y^2}  (x,y_2) \Big| \leq \varepsilon \\
& \qquad \qquad \text{ for all } |y_1|, |y_2| \leq M \text{ with } |y_1-y_2 | \leq \rho.
\end{align*}
\end{assumption}

\begin{assumption} \label{A2.3}
The boundary data satisfies the following regularity property: $g \in L^{q}(\Gamma)$ with $q > \dimension-1.$
\end{assumption}

We observe that the normal derivative $\partial_{\conormal_A} y$ is formally defined by
\[ \partial_{\conormal_A} y = \sum_{i,j=1}^n a_{ij} \partial_{x_i} y(x) \conormal_j(x), \]
where $\conormal(x)$ denotes the outward unit normal vector to $\Gamma$ at the point $x$. Due to the Lipschitz regularity of $\Gamma$ such vector $\conormal(x)$ exists for almost all $x \in \Gamma.$ For a rigorous definition of the normal derivative in a trace sense the reader is referred, for instance, to \cite{Casas-Fernandez1989}.

Throughout this paper the following notation will be used:
\begin{equation} \label{E2.1} \begin{array}{l}
 m_u := \text{ess}\,  \text{inf}_{x \in \Omega} u(x), \quad \Lambda_u := \min\{ \Lambda_A , a_0+m_u\}, \quad \forall u \in L^2(\Omega), \vspace{2mm} \\
\A_0 := \{ u \in L^2(\Omega) : a_0 + m_u >0 \},
\end{array}
\end{equation}
where $\Lambda_A$ and $a_0$ were introduced in Assumptions \ref{A2.1} and \ref{A2.2}, respectively.
It is well known that $H^1(\Omega) \subset L^r(\Omega)$ for every $r \leq \frac{2\dimension}{\dimension-2}$, with $r < \infty$ if $\dimension =2$. Hence, we have,
\begin{equation} \exists C_{r,\Omega} >0  \text{ such that } \|y\|_{L^r(\Omega)} \leq C_{r,\Omega} \|y\|_{H^1(\Omega)}, \quad\forall y \in H^1(\Omega).  \label{E2.2}
\end{equation}
Analogously, since $H^{1/2}(\Gamma)$ is continuously embedded in $L^{q'}(\Gamma)$ for $q > \dimension-1$, where $q' = \frac{q}{q-1}$ denotes the conjugate of $q$, we also have,
\begin{equation} \exists C_{q',\Gamma} >0 \text{ such that } \|y\|_{L^{q'}(\Gamma)} \leq C_{q',\Gamma} \|y\|_{H^1(\Omega)}, \quad\forall y \in H^1(\Omega). \label{E2.3}
\end{equation}

\begin{theorem} \label{T2.1}
There exists $\mu \in (0,1]$ such that for every $u \in \A_0$ there exists a unique solution $y_u \in H^1(\Omega) \cap C^{0,\mu}(\bar \Omega)$ of \eqref{E1.1}. Furthermore, the following estimates hold:
\begin{align}
& \|y_u\|_{H^1(\Omega)} \leq \frac{1}{\Lambda_u} \left ( C_{p',\Omega} \|a(\cdot,0)\|_{L^p(\Omega)} + C_{q',\Gamma} \|g\|_{L^q(\Gamma)} \right ), \label{E2.4} \\
& \|y_u\|_{L^{\infty}(\Omega)} \leq \frac{1}{\Lambda_u} C_{p,q}, \label{E2.5} \\
& \|y_u\|_{C^{0,\mu}(\bar \Omega)} \leq  C_{\mu,\infty} \left ( \|a(\cdot,0)\|_{L^p(\Omega)} + \|u\|_{L^2(\Omega)} + \|g\|_{L^q(\Gamma)} \right ), \label{E2.6}
\end{align}
where $C_{p,q}$ depends on  $\|a(\cdot,0)\|_{L^p(\Omega)}$ and  $\|g\|_{L^q(\Gamma)}$, and $C_{\mu,\infty}$ depends as well on $a(\cdot,0)$ and $g$ and on a monotone nondecreasing way on $\|y_u\|_{L^{\infty}(\Omega)}$.
\end{theorem}
\begin{proof}
We define the mapping
\begin{equation} b: \Omega \times \mathbb R \longrightarrow \mathbb R, \   b(x,y) := a(x,y) - a(x,0) - a_0 y.  \label{E2.7}
\end{equation}
Then, $b$ satisfies $b(y,0)=0$ and $ \frac{\partial b}{\partial y} (x,y) \geq 0$ due to Assumption \ref{A2.2}. Furthermore, \eqref{E1.1} can be written in the form
\begin{equation} \label{E2.8}
\left\{\begin{array}{l} Ay + (a_0+u)y+ b(x,y)= -a(x,0) \ \  \mbox{in } \Omega,\vspace{2mm}\\  \partial_{\conormal_A} y = g\ \ \mbox{on }\Gamma. \end{array}\right.
\end{equation}
From Assumption \ref{A2.1} and \eqref{E2.1}-\eqref{E2.2} we get that  $A+(a_0+u)I : H^1(\Omega) \longrightarrow H^1(\Omega)^*$ is a linear operator satisfying the following properties
\[
\langle (A+(a_0+u))y,y \rangle \geq \Lambda_u \|y\|^2_{H^1(\Omega)} \text{ and } \langle (A+(a_0+u))y,\phi \rangle \leq M_u \|y\|_{H^1(\Omega)} \|\phi \|_{H^1(\Omega)},
\]
where $M_u := M_A + a_0 + C^2_{4,\Omega} \|u\|_{L^2(\Omega)}.$  Therefore, there exists a  unique solution of \eqref{E2.8} in $H^1(\Omega) \cap L^{\infty}(\Omega)$ (see e.g. \cite{Casas93}). Inequality \eqref{E2.4} follows easily by testing \eqref{E2.8} by $y$ and using the established coercivity of the operator $A+(a_0+u)I$ and the fact that $b(x,y)y \geq 0$ together with \eqref{E2.2} and \eqref{E2.3}. To deduce \eqref{E2.5} similar to \cite{Alibert-Raymond97} we introduce the function $y_k(x,t):=y(x,t)-\proj_{[-k,k]}(y(x,t))$  for every integer $k \geq 1.$
Testing \eqref{E2.8} with $y_k$, using that $\partial_{x_i}y \partial_{x_j}y_k = \partial_{x_i} y_k \partial_{x_j} y_k$, the inequality $(a_0+a)yy_k \geq (a_0+a)y^2_k$, and $b(x,y)y_k \geq 0$ we infer,
\[ \|y_k \|^2_{H^1(\Omega)} \leq \frac{1}{\Lambda_u} \left ( \int_\Omega | a(\cdot,0) | |y_k(x) | \dx + \int_\Gamma |g(x) | |y_k(x) | \dx \right ). \]
Following the techniques of \cite{Stampacchia65} and \cite{Alibert-Raymond97} we deduce the estimate. For the last estimate,  we write \eqref{E1.1} in the form
\begin{equation*}
\left\{\begin{array}{l} Ay + y= (1-u)y - a(x,y) \ \  \mbox{in } \Omega,\vspace{2mm}\\  \partial_{\conormal_A} y = g\ \ \mbox{on }\Gamma. \end{array}\right.
\end{equation*}
and we denote by $M:= \|y\|_{L^{\infty}(\Omega)}.$ With Assumption \ref{A2.2} and the mean value Theorem we deduce,
\[ |a(x,y) | \leq | a(x,0) | + C_{a,M}M. \]
In addition, we have $\|(1-u) y\|_{L^2(\Omega)} \leq ( \|u\|_{L^2(\Omega)} + \sqrt{|\Omega|} ) M.$ Combining these estimates with the results of \cite{Nittka11} the existence of $\mu \in (0,1]$ such that the regularity $y \in C^{0,\mu}(\bar \Omega)$ and inequality \eqref{E2.6} follow. \qed
\end{proof}

Next we consider the diffentiability of the mapping $u \to y_u$.
\begin{theorem} \label{T2.2}
There exists an open set $\A$ in $L^2(\Omega)$ such that $\A_0 \subset \A$ and $\forall u \in \A$ the equation \eqref{E1.1} has a unique solution $y_u \in H^1(\Omega) \cap C^{0,\mu}(\bar \Omega)$, where $\mu \in (0,1]$ was introduced in Theorem \ref{T2.1}.
Further, the mapping $G : \A \longrightarrow H^1(\Omega) \cap C^{0,\mu} (\bar \Omega)$ defined by $G(u):=y_u$ is of class $C^2$ and $\forall u \in \A$ and $\forall v, v_1, v_2 \in L^2(\Omega)$ the functions
$z =G'(u)v$ and $w=G''(u)(v_1,v_2)$ are the unique solutions of the equations:
\begin{align}
& \left\{\begin{array}{l} Az + \displaystyle \frac{\partial a}{\partial y}(x,y_u)z + uz + vy_u= 0\ \  \mbox{in } \Omega,\vspace{2mm}\\  \partial_{\conormal_A} z = 0\ \ \mbox{on }\Gamma, \end{array}\right.
\label{E2.9} \\
& \left\{\begin{array}{l} Aw + \displaystyle \frac{\partial a}{\partial y}(x,y_u)w + \displaystyle \frac{\partial^2 a}{\partial y^2}(x,y_u)z_{u,v_1}z_{u,v_2} + uw + v_1z_{u,v_2} + v_2z_{u,v_1} = 0\ \  \mbox{in } \Omega,\vspace{2mm}\\  \partial_{\conormal_A} w = 0\ \ \mbox{on }\Gamma, \end{array}\right.
\label{E2.10}
\end{align}
where $z_{u,v_i} = G'(u)v_i$, $i=1,2$.
\end{theorem}
\begin{proof}
We define the space
\[ Y_{A} := \{ y \in H^1(\Omega) \cap C^{0,\mu}(\bar \Omega) : Ay \in L^p(\Omega) , \  \partial_{\conormal_A} y \in L^q(\Gamma) \}  \]
which is a Banach space when endowed with the graph norm.
We also define the mapping
\[ {\mathcal F}: L^2(\Omega) \times Y_A \longrightarrow L^p(\Omega) \times L^q(\Gamma), \quad {\mathcal F}(u,y) := ( Ay + a(x,y) + uy, \partial_{\conormal_A} y - g ) \]
From Assumption \ref{A2.2} we deduce that $\mathcal F$ is of class $C^2$. For every $(\bar u,\bar y)\in\A_0\times Y_A$ the derivative $\frac{\partial {\mathcal F}}{\partial y}(\bar u,\bar y) : Y_A \longrightarrow L^p(\Omega) \times L^q(\Gamma)$, given by
\[ \frac{\partial {\mathcal F}}{\partial y}(\bar u,\bar y)z = \left  ( Az + \displaystyle \frac{\partial a}{\partial y}(x,\bar y)z + \bar u z ,    \partial_{\conormal_A} z \right ) \ \forall z \in Y_{A},\]
is linear and continuous.
Using Theorem \ref{T2.1}, we deduce that the equation
 \begin{equation} \label{E2.11}
 \left\{\begin{array}{l} Az + \displaystyle \frac{\partial a}{\partial y}(x,\bar y)z + \bar u z= f\ \  \mbox{in } \Omega,\vspace{2mm}\\  \partial_{\conormal_A} z = h\ \ \mbox{on }\Gamma, \end{array}\right.
\end{equation}
has unique solution $z \in Y_A$ for all $(f,h) \in L^p(\Omega) \times L^q(\Gamma)$.
The open mapping theorem implies that  $\frac{\partial {\mathcal F}}{\partial y}(\bar u,\bar y)$ is an isomorphism and
%
there exists $ \varepsilon_{\bar u} >0$ and $\varepsilon_{\bar y} >0$, such that $\forall u \in B_{\varepsilon_{\bar u}}(\bar u) \subset L^2(\Omega)$ the equation ${\mathcal F}(u,y) = 0$ has a unique solution $y_u$ in the ball $B_{\varepsilon_{\bar y}}(\bar y) \subset Y.$  Moreover the mapping $u \in  B_{\varepsilon_{\bar u}}(\bar u) \to y_u \in B_{\varepsilon_{\bar y}}(\bar y)$ is of class $C^2$. Without loss of generality, we assume $\varepsilon_{\bar u} < \frac{\Lambda_{\bar u}}{C^2_{4,\Omega}}$, $\Lambda_{\bar u}$ was defined in \eqref{E2.1} and was $C_{4,\Omega}$ introduced in \eqref{E2.2} for $r=4.$
We prove that for every $u \in B_{\varepsilon_{\bar u}}$ the equation ${\mathcal F}(u,y) = 0$ has unique solution $y \in Y_A.$ Indeed, suppose that $y_1,y_2$ are two solutions of
${\mathcal F}(u,y) = 0.$ We set $y=y_1-y_2$, subtract the corresponding equations, and apply the mean value theorem to deduce that $y$ satisfies,
 \begin{equation} \label{E2.12}
 \left\{\begin{array}{l} Ay + \displaystyle \frac{\partial a}{\partial y}(x,y_1+\theta_x y) y + u y= 0\ \  \mbox{in } \Omega,\vspace{2mm}\\  \partial_{\conormal_A} y = 0\ \ \mbox{on }\Gamma, \end{array}\right.
\end{equation}
where $\theta_x:\Omega\to[0,1]$ is a measurable function.
The equation \eqref{E2.12} can be written as
 \begin{equation} \label{E2.13}
 \left\{\begin{array}{l} Ay + \displaystyle [ \frac{\partial a}{\partial y}(x,y_1+\theta_x y) + \bar u ] y + (u-\bar u) y= 0\ \  \mbox{in } \Omega,\vspace{2mm}\\  \partial_{\conormal_A} y = 0\ \ \mbox{on }\Gamma. \end{array}\right.
\end{equation}
Testing \eqref{E2.13} with $y$ we get,
\[ \left ( \Lambda_{\bar u} - C^2_{4,\Omega} \varepsilon_{\bar u} \right ) \|y\|^2_{H^1(\Omega)} \leq  \Lambda_{\bar u} \|y\|^2_{H^1(\Omega)} - C^2_{4,\Omega} \|u-\bar u\|_{L^2(\Omega)} \|y\|^2_{H^1(\Omega)} \leq 0. \]
Hence, $y=0$ holds. Finally, defining in $L^2(\Omega)$ the open set $\A = \cup_{\bar u \in \A_0} B_{\varepsilon_{\bar u}} (\bar u)$ and the mapping $G: \A \longrightarrow Y$ such that $G(u)=y_u,$ we have that $G$ is of class of $C^2$. Moreover, the equations \eqref{E2.9} and \eqref{E2.10} are obtained differentiating with respect to $u$ the identity ${\mathcal F}(u,G(u)) = 0$. \qed
\end{proof}

\section{Analysis of the optimal control problem} \label{S3}
\setcounter{equation}{0}
In this Section we are going to prove existence of solutions of problem $\Pb$ and we analyze the first and second order optimality conditions for a local minimizer.  Along this paper a local minimizer is understood in $L^2(\Omega)$ sense. To this end we make the following hypothesis.
\begin{assumption} \label{A3.1} We assume that the conditions $\tichonov >0$ and $ - a_0 < \umin < \umax \leq \infty$ hold.
\end{assumption}
\begin{assumption} \label{A3.2} The function $L: \Omega \times \mathbb R \longrightarrow \mathbb R$ is Carath\'eodory and of class of $C^2$ with respect to the second variable.
Further the following properties hold for almost all $x \in \Omega$:
\begin{align*}
& \bullet L(\cdot, 0) \in L^1(\Omega) \\
& \bullet \forall M >0, \  \exists L_M \in L^p(\Omega) \text{ such that } \Big|  \frac{\partial L}{\partial y}(x,y) \Big| \leq L_M(x) \ \forall |y| \leq M, \\
& \bullet \forall M >0, \ \exists C_{L,M} \in \mathbb R   \text{ such that } \Big|  \frac{\partial^2 L}{\partial y^2}(x,y) \Big| \leq C_{L,M} \ \forall |y| \leq M, \\
& \bullet \forall \varepsilon >0  \text{ and } \forall M >0 \ \exists \rho >0 \text{ such that } \\
& \qquad\qquad  \Big|  \frac{\partial^2 L}{\partial y^2}(x,y_1) - \frac{\partial^2 L}{\partial y^2}(x,y_2)  \Big| \leq \varepsilon \ \forall |y_1|, |y_2| \leq M \text{ with } |y_1-y_2| \leq \rho.
\end{align*}
\end{assumption}
\begin{remark} \label{R3.1}
1. From Assumption \ref{A3.1} and \eqref{E2.1} we deduce that $a_0 + m_u \geq a_0 + \umin >0 \ \forall u \in \uad.$ Consequently, we have that $\uad \subset \A_0 \subset \A$. \vspace{2mm}\\
2. Defining $\Lambda = \min\{ \Lambda_A,a_0+ \umin \}$ we infer that $0 < \Lambda \leq \Lambda_u$ for every $u \in \uad.$ Consequently, \eqref{E2.4} and \eqref{E2.5} hold with $\Lambda$ instead of $\Lambda_u$, and $C_{\mu,\infty}$ in \eqref{E2.6} can be chosen independently of $u \in \uad$.
 \end{remark}
As a consequence of Theorem \ref{T2.2} and Assumption \ref{A3.2} we deduce the differentiability of functional $J$.
\begin{theorem} \label{T3.1} The functional $J: \A \longrightarrow \mathbb R$ is of class $C^2$ and its derivatives are given by the expressions:
\begin{align}
& J'(u)v = \int_\Omega (\tichonov u - y_u \varphi_u)v \dx \ \forall u \in \A , \ \forall v \in L^2(\Omega), \label{E3.1} \\
& J''(u)(v_1,v_2) = \int_\Omega \Big[ \frac{\partial^2 L}{\partial y^2}(x,y_u) - \varphi_u  \frac{\partial^2 a}{\partial y^2}(x,y_u) \Big] z_{u,v_1} z_{u,v_2}\dx \label{E3.2} \\
& \qquad - \int_\Omega \Big[ v_1z_{u,v_2} + v_2 z_{u,v_1} \Big] \varphi_u \dx + \tichonov \int_\Omega v_1 v_2 \dx, \ \forall u \in \A , \ \forall v_1, v_2 \in L^2(\Omega), \nonumber
\end{align}
where $z_{u,v_i} = G'(u)v_i$, $i=1,2$ and $\varphi_u \in H^1(\Omega) \cap C^{0,\mu}(\bar \Omega)$ is the adjoint state, the unique solution of the equation
\begin{equation}
 \left\{\begin{array}{l} A^*\varphi + \displaystyle \frac{\partial a}{\partial y}(x,y_u)\varphi + u\varphi = \frac{\partial L}{\partial y}(x,y_u) \ \  \mbox{in } \Omega,\vspace{2mm}\\  \partial_{\conormal_{A^*}} \varphi = 0\ \ \mbox{on }\Gamma. \end{array}\right.
\label{E3.3}
\end{equation}
\end{theorem}
\begin{proof}
First let us analyze \eqref{E3.3}. To prove existence, uniqueness and regularity of solution \eqref{E3.3} we first observe that there exists $\bar u \in \A_0$ such that $u \in B_{\varepsilon_{\bar u}}(\bar u)$, where $\varepsilon_{\bar u}$ is defined in the proof of Theorem \ref{T2.2}. Then, \eqref{E3.3} can be written as
\begin{equation*}
 \left\{\begin{array}{l} A^*\varphi + \displaystyle \Big[ \frac{\partial a}{\partial y}(x,y_u) + \bar u \Big] \varphi + (u-\bar u) \varphi = \frac{\partial L}{\partial y}(x,y_u) \ \  \mbox{in } \Omega,\vspace{2mm}\\  \partial_{\conormal_{A^*}} \varphi = 0\ \ \mbox{on }\Gamma. \end{array}\right.
\end{equation*}
Setting $M= \|y_u\|_{L^{\infty}(\Omega)}$, from Assumption \ref{A3.2} we obtain $ \Big| \frac{\partial L}{\partial y}(x,y_u)\Big| \leq L_M(x)$ with $L_M \in L^p(\Omega)$. Then, arguing  as for \eqref{E2.13} we obtain the coercivity of the linear equation. The existence and uniqueness of a solution in $H^1(\Omega)$ follows from Lax-Milgram theorem. Finally, using again \cite{Nittka11} we deduce the $C^{0,\mu}(\bar \Omega)$ regularity.

The fact that $J$ is of class $C^2$ is an immediate consequence of the chain rule, Theorem \ref{T2.2}, and Assumption \ref{A3.2}. Moreover, we have
\begin{align*}
& J'(u)v= \int_\Omega \Big[\frac{\partial L}{\partial y}(x,y_u) z_{u,v} + \tichonov uv \Big] \dx, \\
& J''(u)(v_1,v_2) = \int_\Omega \Big[ \frac{\partial L}{\partial y}(x,y_u) w +  \frac{\partial^2 L}{\partial y^2}(x,y_u)z_{u,v_1} z_{u,v_2} + \tichonov v_1v_2 \Big] \dx,
\end{align*}
where $z_{u,v} = G'(u)v$, $z_{u,v_i} = G'(u)v_i$, $i=1,2$, and $w=G''(u)(v_1,v_2).$ Combining these expressions with \eqref{E2.9}, \eqref{E2.10}, and \eqref{E3.3} the formulas \eqref{E3.1} and \eqref{E3.2} follow.
\qed
\end{proof}
In the above theorem we have proved that the mapping $\Phi : \A \longrightarrow H^1(\Omega) \cap C^{0,\mu}(\bar \Omega)$ given by $\Phi (u):= \varphi_u$ is well defined. In the next theorem its differentiability is established.
\begin{theorem} \label{T3.2}
The mapping $\Phi$ is of class $C^1$ and for all $u\in \A$ and $v\in L^2(\Omega)$ the function $\eta_{u,v} = \Phi'(u)v$ is the unique solution of
\begin{equation}
 \left\{\begin{array}{l} A^*\eta + \displaystyle \frac{\partial a}{\partial y}(x,y_u)\eta + u\eta = \Big[\frac{\partial^2 L}{\partial y^2}(x,y_u)-\varphi_u  \frac{\partial^2 a}{\partial y^2}(x,y_u)\Big] z_{u,v} - v\varphi_u \ \  \mbox{in } \Omega,\vspace{2mm}\\  \partial_{\conormal_{A^*}} \eta = 0\ \ \mbox{on }\Gamma, \end{array}\right.
\label{E3.4}
\end{equation}
where $z_{u,v}= G'(u)v.$
\end{theorem}
\begin{proof}
According to Assumption \ref{A3.2} and the fact that $y_u,\varphi_u,z_{u,v} \in L^{\infty}(\Omega)$ we deduce that the right hand side of \eqref{E3.4} belongs to $L^2(\Omega)$. As for \eqref{E3.3} the existence, uniqueness, and regularity of $\eta_{u,v}$ follows. To prove the differentiability of $\Phi$ we define
\[Y_{A^*} = \{ \varphi \in H^1(\Omega) \cap C^{0,\mu}(\bar \Omega) : A^* \varphi \in L^p(\Omega) \text{ and } \partial_{\conormal_{A^*}} \varphi = 0 \}\] and
$\mathcal G : \A \times Y_{A^*} \longrightarrow L^p(\Omega)$ by
\[ {\mathcal G}(u,\varphi ) :=A^* \varphi + \frac{\partial a}{\partial y}(x,y_u) \varphi + u\varphi - \frac{\partial L}{\partial y}(x,y_u). \]
From Assumptions \ref{A2.2} and \ref{A3.2} we deduce that $\mathcal G$ is a $C^1$ mapping. We have  that
\[
\frac{\partial \mathcal G}{\partial \varphi}(u,\varphi) \eta = A^* \eta + \frac{\partial a}{\partial y}(x,y_u) \eta +u \eta
\]
and $\frac{\partial \mathcal G}{\partial \varphi}(u,\varphi) : Y_{A^*} \longrightarrow L^p(\Omega)$ is an isomorphism. Then, applying the implicit function theorem and differentiating the identity $\mathcal{G}(u,\Phi(u)) = 0$ the result follows.  \qed
\end{proof}

The following corollary is a straightforward application of formula \eqref{E3.2} and equation \eqref{E3.4}.
\begin{corollary} \label{C3.1} For every $v_1, v_2 \in L^2(\Omega)$ and all $u \in \A$, the following identities hold
\begin{align*}
J''(u)(v_1,v_2) & = \int_{\Omega}\Big[ \tichonov v_1 - (\varphi_u z_{u,v_1}  + y_u \eta_{u,v_1}) \Big] v_2 \dx \\
&=  \int_{\Omega}\Big[ \tichonov v_2 - (\varphi_u z_{u,v_2}  + y_u \eta_{u,v_2}) \Big] v_1 \dx
\end{align*}
\end{corollary}

\begin{theorem} \label{T3.3} Problem $\Pb$ has at least one solution. Moreover, if $\bar u \in \uad$ is a local minimizer of $\Pb$ then there exist $\bar y, \bar \varphi \in H^1(\Omega) \cap C^{0,\mu} (\bar \Omega)$ such that
\begin{align}
& \left\{\begin{array}{l} A\bar y + a(x,\bar y) + \bar u \bar y= 0\ \  \mbox{in } \Omega,\vspace{2mm}\\  \partial_{\conormal_A} \bar y = g\ \ \mbox{on }\Gamma. \end{array}\right.
\label{E3.5} \\
& \left\{\begin{array}{l} A^* \bar \varphi + \displaystyle \frac{\partial a}{\partial y}(x,\bar y)\bar \varphi + \bar u \bar \varphi = \frac{\partial L}{\partial y}(x,\bar y) \ \  \mbox{in } \Omega,\vspace{2mm}\\  \partial_{\conormal_{A^*}} \bar  \varphi = 0\ \ \mbox{on }\Gamma. \end{array}\right. \label{E3.6} \\
& \bar u(x) = \proj_{[\umin, \umax]} \left ( \frac{1}{\tichonov} \bar y(x) \bar \varphi(x) \right ). \label{E3.7}
\end{align}
\end{theorem}
The existence of a solution follows by usual arguments, taking a minimizing sequence, and observing that if $u_k \rightharpoonup \bar u$ in $L^2(\Omega)$ then $y_{u_k} \to \bar y = y_{\bar u}$ strongly in $H^1(\Omega) \cap C(\bar \Omega)$. This statement is an immediate consequence of estimates \eqref{E2.4} - \eqref{E2.6} and Remark \ref{R3.1}. The optimality system follows from \eqref{E3.1}, \eqref{E3.3}, and the fact that $\uad$ is convex.

From now on $(\bar u,\bar y,\bar \varphi) \in \uad \times [ H^1(\Omega) \cap C^{0,\mu}(\bar \Omega) ]^2$ will denote a triplet that satisfies \eqref{E3.5}-\eqref{E3.7}.
Associated with this triplet we define the cone of critical directions
\begin{equation} \label{E3.8}
C_{\bar u} = \{ v \in L^2(\Omega) : v(x) {=} 0 \text{ if } \tichonov \bar u(x) - \bar y(x) \bar \varphi(x) {\ne} 0  \text{ a.e. in $\Omega$ and (\ref{E3.9})  holds} \},
\end{equation}
\begin{equation} \label{E3.9}
v(x)  \left\{ \begin{array}{ll} \geq 0 & \text{ if } \bar u(x) = \umin, \\ \leq 0 & \text{ if } \bar u(x) = \umax. \end{array} \right.
\end{equation}

Regarding the second order optimality conditions we have the following result.
\begin{theorem} \label{T3.4}
If $\bar u$ is a local minimizer of $\Pb$, then $J''(\bar u)v^2 \geq 0 \ \forall v \in C_{\bar u}$ holds. Conversely, if $\bar u \in \uad$ satisfies the first order optimality conditions \eqref{E3.5}--\eqref{E3.7} and $J''(\bar u)v^2 > 0 \ \forall v \in C_{\bar u} \setminus\{0\}$, then there exist $\varepsilon >0$ and $\delta >0$ such that
\begin{equation} \label{E3.10}
J(\bar u) + \frac{\delta}{2} \|u - \bar u\|^2_{L^2(\Omega)} \leq J(u) \quad\forall u \in \uad \text{ with } \|u-\bar u\|_{L^2(\Omega)} \leq \varepsilon.
\end{equation}
\end{theorem}

The proof of this theorem is standard. The reader is referred, for instance, to \cite{Casas-Troltzsch12}. For the subsequent analysis the strict complementarity condition will be needed.
\begin{definition} \label{D3.1}
Let us define
\[ \Sigma_{\bar u} = \{ x \in \Omega : \bar u(x) \in \{ \umin, \umax \} \text{ and } \tichonov \bar u(x) - \bar y(x) \bar \varphi (x) = 0 \}. \]
We say that the strict complementarity condition is satisfied at $\bar u$ if $ | \Sigma_{\bar u} | = 0,$ where $ | \cdot | $ stands for the Lebesgue measure.
\end{definition}
This notion is an extension to the case of infinite constraints of the usual strict complementarity condition in finite dimensional nonlinear programming.

For every $\tau \geq 0$, we define the subspace
\begin{equation} \label{E3.11}
 T^{\tau}_{\bar u} = \{ v \in L^2(\Omega) : v(x) = 0 \text{ if } | \tichonov \bar u(x) - \bar y (x) \varphi (x) | > \tau \text{ a.e. in $\Omega$} \}.
\end{equation}
If $\tau =0$ we simply denote $ T_{\bar u} = T^{0}_{\bar u}.$
\begin{theorem} \label{T3.5}
Assume that $\bar u$ satisfies the strict complementarity condition. Then, the following properties hold: \\
1- $T_{\bar u} = C_{\bar u}.$ \\
2- If $\bar u$ satisfies the second order optimality condition $J''(\bar u)v^2 >0 \ \forall v \in C_{\bar u} \setminus\{0\}$, then
\begin{equation} \label{E3.12}
\exists \tau >0 \text{ and }  \kappa >0 \text{ such that } J''(\bar u) v^2 \geq \kappa \|v\|^2_{L^2(\Omega)} \ \forall v \in T^\tau_{\bar u}.
\end{equation}
\end{theorem}
\begin{proof}
1- It is obvious that $C_{\bar u} \subset T_{\bar u}$. Let us prove the converse inclusion. If $v \in T_{\bar u}$ we have to prove that $v$ satisfies the sign conditions \eqref{E3.9}. If $\bar u(x) = \umin,$  then from \eqref{E3.7} we deduce $\tichonov \bar u(x) - \bar y(x) \bar \varphi(x) \geq 0.$ Hence, with the strict complementarity condition we get that
$\tichonov \bar u(x) - \bar y(x) \bar \varphi(x) > 0$  for almost all  $x$ such that $\bar u(x) = \umin$. Since $v \in T_{\bar u}$ we conclude $v(x)=0$ for almost all $x$ such that $\bar u(x) = \umin$. In a similar way we argue when $\bar u(x) = \umax.$ \vspace{2mm} \\
2- We argue by contradiction. If the statement is false, then $\forall k \geq 1$ $\exists v_k \in T^{1/k}_{\bar u}$ such that $J''(\bar u)v^2_k < \frac{1}{k} \|v_k\|^2_{L^2(\Omega)}$.
Dividing $v_k$ by $\|v_k\|_{L^2(\Omega)}$ and denoting the result again by $v_k$, we obtain
\begin{equation} \label{E3.13} v_k \in  T^{1/k}_{\bar u}, \ \|v_k\|_{L^2(\Omega)} = 1, \text{ and } J''(\bar u)v^2_k < \frac{1}{k}. \end{equation}
Taking a subsequence that we denote again by $v_k$ we have that $v_k \rightharpoonup v$ in $L^2(\Omega)$.  The rest of the proof is divided in three steps. \vspace{2mm} \\
{\it Step} 1- $v \in C_{\bar u}$. According to statement 1 of the theorem we only need to prove that $v(x) =0$ if $\tichonov \bar u(x)-\bar y(x)\bar\varphi(x) \ne 0.$ Given $\varepsilon >0$ we define
$\Omega_\varepsilon =\{ x \in \Omega : | \tichonov \bar u(x)-\bar y(x) \bar \varphi (x) | > \varepsilon \}. $  By definition of $T^{1/k}_{\bar u},$ we have that $v_k(x) = 0$ a.e. in $\Omega_\varepsilon$ $ \forall k > \frac{1}{\varepsilon}$. Therefore, $v(x) = 0$ a.e in $\Omega_\varepsilon$ holds. Since $\varepsilon>0$ can be selected arbitrarily small we conclude that $v(x) = 0$ for almost all $x$ in $\Omega$ such that $\tichonov \bar u(x)-\bar y(x) \bar \varphi(x) \ne 0,$ hence $v \in T_{\bar u} = C_{\bar u}$. \vspace{2mm} \\
{\it Step} 2- $J''(\bar u)v^2 \leq 0.$ Since $v_k \rightharpoonup v$ in $L^2(\Omega)$ we get that $z_{\bar u,v_k} = G'(\bar u)v_k \rightharpoonup G'(\bar u)v = z_{\bar u,v}$ in $H^1(\Omega) \cap C^{0,\mu} (\bar \Omega).$ Therefore, the convergence $z_{\bar u,v_k} \to z_{\bar u,v}$ is strong in $C(\bar \Omega).$ Then, we easily pass to the limit in \eqref{E3.2} to deduce with \eqref{E3.13} that $J''(\bar u)v^2 \leq \lim \inf _{k \to \infty} J''(\bar u) v^2_k \leq 0. $ \vspace{2mm} \\
{\it Step 3- Final Contradiction.} From Steps 1 and 2, and the fact that $\bar u$ satisfies the second order sufficient optimality condition, we deduce that $v=0$ and consequently, $z_{\bar u,v_k} \to 0 \text{ in } C(\bar \Omega).$  Using the fact that $\|v_k\|_{L^2(\Omega)} = 1,$ we have that  $J''(\bar u)v^2_k = \tichonov + \varepsilon_k$ with $\varepsilon_k \to 0$. Therefore,
$\tichonov = \lim_{k \to \infty} J''(\bar u)v^2_k \leq 0$ which contradicts the strict positivity of $\tichonov$. \qed
\end{proof}

\section{Convergence of the semismooth Newton method} \label{S4}
\setcounter{equation}{0}
Following \cite[Chapter 3]{Ulbrich11} we are going to describe the abstract framework where our numerical algorithm fits.

\begin{definition} \label{D4.1}
Given two Banach spaces $U$ and $X$, an open subset $\A$ of $U$, a continuous function $F:\A \longrightarrow X$, and a set-value mapping $\partial F : \A \longrightarrow {\mathcal L}(U,X)$ such that $\partial F(u) \neq \emptyset \ \forall u \in \A$, we say that $F$ is $\partial F$-semismooth at $\bar u \in \A$ if
\[ \lim_{v \to 0} \sup_{M \in \partial F(\bar u + v)} \frac{ \|F(\bar u + v) - F(\bar u) - Mv ||_{X}}{ \|v\|_{U} } = 0. \]
$F$ is said $\partial F$-semismooth at $\A$ if it is $\partial F$-semismooth at every $u \in \A.$
\end{definition}

The abstract semismooth Newton algorithm takes the form:

\LinesNumbered
\begin{algorithm2e}[h!]
\caption{Semismooth Newton method.}\label{Alg1}
\DontPrintSemicolon
Initialize Choose $u_0\in \A$. Set $j=0$.\;
\For{$j\geq 0$}{
Choose $M_j\in\partial F(u_j)$ and solve $M_jv_j=-F(u_j)$\label{line3}.\;
Set $u_{j+1}=u_j+v_j$ and $j=j+1$.\;
}
\end{algorithm2e}

\begin{theorem} \label{T4.1}
Suppose that $F: \A \longrightarrow X$ is $\partial F$-semismooth at $\bar u \in \A$ solution of $F(u) = 0$ locally unique. Assume, furthermore, that $\forall j$ the operator $M_j \in \partial F(u_j)$ is an isomorthism and there exists $C_{F}>0$ such that
\begin{equation} \label{E4.1} \|M^{-1}_j\|_{{\mathcal L}(X,U)} \leq C_F \quad \forall j \geq 0. \end{equation}
Then, there exists $\delta >0$ such that for all $u_0 \in \A$ with $\|u_0-\bar u\|_{U} \leq \delta$ the sequence $\{u_j\}_{j \geq 0}$ generated by the semismooth Newton method (Algorithm $1$) converges superlinearly to $\bar u$.
\end{theorem}
The proof is this theorem can be found in \cite[Theorem 3.13]{Ulbrich11}; see also \cite{Hintermuller-Ito-Kunisch03}.
Let us put our problem into this particular framework. Let $X=U=L^2(\Omega)$ and $\A$ be the open set introduced in Theorem \ref{T2.2}. We define by $F: \A \longrightarrow L^2(\Omega)$  by
\[ F(u) = u - \proj_{[\umin,\umax]} \left (\frac{1}{\tichonov} y_u \varphi_u \right ). \]
Due to Theorem \ref{T3.3} any local minimizer of $\Pb$ is a solution of $F(u) = 0.$  In order to define $\partial F(u) \ \forall u \in \A$ we introduce some additional functions.
\begin{align*}
& S : \A  \longrightarrow L^2(\Omega), \quad S(u)= \frac{1}{\tichonov} G(u) \Phi (u), \\
& \psi : \mathbb R  \longrightarrow \mathbb R, \quad  \psi(t) = \proj_{[\umin,\umax]} (t), \\
& \Psi :  \A \longrightarrow  L^2(\Omega), \quad  \Psi(u)(x) = \psi(S(u)(x)).
\end{align*}
For every $u \in \A$ we define
\begin{align*} \partial \Psi (u) & =  \big \{ N \in {\mathcal L}(L^2(\Omega),L^2(\Omega)) : Nv = h S'(u)v \ \forall v \in L^2(\Omega) \text{ and for some measurable } \\
& \qquad \text{ function } h: \Omega \longrightarrow \mathbb R \text{ such that }  h(x) \in \partial \psi (S(u)(x))\big \}.
\end{align*}
We observe that $\psi$ is a Lipschitz function and by $\partial \psi (t)$ we denote the subdifferential in Clarke's sense; see \cite[Chapter 2]{Clarke83}. Note that
\[ \partial \psi(t) = \left\{ \begin{array}{cl} \{1\} & \text{ if } t \in (\umin,\umax), \\ \{0\} & \text{ if } t \not\in [\umin, \umax], \\  {[0,1]} & \text{ if } t \in \{\umin,\umax \}.
 \end{array}\right. \]
According to \cite[Proposition 2.26]{Ulbrich11}, $\psi$ is $1$-order $\partial \psi$-semismooth.

\begin{theorem} \label{T4.3}
$\Psi$ is $\partial \Psi$-semismooth in $\A$.
\end{theorem}

\begin{proof}
Since $\Psi$ is a superposition operator of $\psi$ and $S$, we will apply \cite[Theorem 3.49]{Ulbrich11} to deduce that $\partial \Psi$-semismooth in $\A$.
To this end it is enough to prove that $ S : \A  \longrightarrow L^2(\Omega)$ is $C^1$ and that $ S : \A  \longrightarrow L^6(\Omega)$ is locally Lipschitz. The first condition is an immediate consequence of Theorems \ref{T2.2} and \ref{T3.2}. Indeed, since $S(u) = \frac{1}{\tichonov} G(u) \Phi(u)$ we have that
\[S'(u)= \frac{1}{\tichonov} [ G'(u)v \Phi(u) + G(u) \Phi'(u)v ] = \frac{1}{\tichonov} [ z_{u,v}\varphi_u + y_u \eta_{u,v} ]. \]
Let us prove the second condition. Given $u \in \A$, by definition of $\A$ there exists $B_{\varepsilon_{\bar u}}(\bar u) \subset \A$ with $\bar u \in \A_0,$ such that $u \in B_{\varepsilon_{\bar u}}(\bar u)$.
We recall that $\varepsilon_{\bar u}$ was selected such that $C^2_{4,\Omega} \varepsilon_{\bar u} < \Lambda_{\bar u}$; see the proof of Theorem \ref{T2.2}. As a first step we are going to prove that there exists $M>0$ such that $\|y_u\|_{H^1(\Omega) \cap C^{0,\mu}(\bar\Omega)} \leq M \ \forall u \in B_{\varepsilon_{\bar u}}(\bar u)$. Given $u \in B_{\varepsilon_{\bar u}}(\bar u)$, we define $y=y_{\bar u}- y_u$. Subtracting the equations satisfied by $y_{\bar u}$ and $y_u$ and the using mean value theorem we obtain
\[\left\{\begin{array}{l} Ay + \displaystyle \frac{\partial a}{\partial y}(x,y_\theta)y + u y= (u-\bar u) y_{\bar u}\ \  \mbox{in } \Omega,\vspace{2mm}\\  \partial_{\conormal_A} y= 0\ \ \mbox{on }\Gamma, \end{array}\right.\]
where $y_\theta = y_u + \theta y$ for some measurable function $\theta : \Omega \longrightarrow [0,1]$. We rewrite the above equation in the form
\[
Ay + \Big[\frac{\partial a}{\partial y}(x,y_\theta) + \bar u\Big]y + (u - \bar u)y= (u-\bar u) y_{\bar u}.
\]
Testing this equation with $y$ and using that $u \in B_{\varepsilon_{\bar u}}(\bar u)$ we infer
\[
(\Lambda_{\bar u} - C^2_{4,\Omega}\varepsilon_{\bar u})\|y\|_{H^1(\Omega)} \le \varepsilon\|y_{\bar u}\|_{L^2(\Omega)}.
\]
Therefore, we have that $\|y_u\|_{H^1(\Omega)} \le \|y_{\bar u}\|_{H^1(\Omega)} + \|y\|_{H^1(\Omega)} \le C$ for all $u \in B_{\varepsilon_{\bar u}}(\bar u)$ and some real constant $C$.
Now, applying again \cite{Nittka11} to the equation
\[
Ay + \Big[\frac{\partial a}{\partial y}(x,y_\theta) + \bar u\Big]y = (u-\bar u) (y_{\bar u} - y)
\]
we infer that $\|y_u\|_{H^1(\Omega) \cap C^{0,\mu}(\bar\Omega)} \le \|\bar y_u\|_{H^1(\Omega) \cap C^{0,\mu}(\bar\Omega)} + \|y\|_{H^1(\Omega) \cap C^{0,\mu}(\bar\Omega)} \le M$ for all $u \in B_{\varepsilon_{\bar u}}(\bar u)$.

Now, we prove that $G:\A \longrightarrow H^1(\Omega) \cap C^{0,\mu}(\bar\Omega)$ is Lipschitz in the ball $B_{\varepsilon_{\bar u}}(\bar u)$. Of course, we also have the same Lipschitz property from $B_{\varepsilon_{\bar u}}(\bar u)$ to $L^6(\Omega)$. Given $u_1, u_2 \in B_{\varepsilon_{\bar u}}(\bar u)$, we deduce from the generalized mean value theorem
\[
\|G(u_2) - G(u_1)\|_{H^1(\Omega) \cap C^{0,\mu}(\bar\Omega)} \le \sup_{\rho \in [0,1]}\|G'(u_1 + \rho(u_2 - u_1))\|\|u_2 - u_1\|_{L^2(\Omega)},
\]
where $\|G'(u)\|$ stands for the norm in the space $\mathcal{L}(L^2(\Omega),H^1(\Omega)\cap C^{0,\mu}(\bar\Omega))$. Since $u_1 + \rho(u_2 - u_1) \in B_{\varepsilon_{\bar u}}(\bar u)$ for all $\rho \in [0,1]$, the above inequality implies the Lipschitz property if we prove the existence of a constant $C$ such that $\|G'(u)\| \le C$ for every $u \in B_{\varepsilon_{\bar u}}(\bar u)$. We observe that $\|G'(u)\| = \sup_{\|v\|_{L^2(\Omega)} \le 1}\|z_{u,v}\|_{H^1(\Omega)\cap C^{0,\mu}(\bar\Omega)}$. Rewriting the equation satisfied by $z_{u,v}$ in the way
\[
Az_{u,v} + \Big[\frac{\partial a}{\partial y}(x,y_u) + \bar u\Big]z_{u,v} + (u - \bar u)z_{u,v} = -vy_u,
\]
using that $\|y_u\|_{C(\bar\Omega)} \le M$, and testing the above equation with $z_{u,v}$ we infer
\[
(\Lambda_{\bar u} - C^2_{4,\Omega}\varepsilon_{\bar u})\|z_{u,v}\|^2_{H^1(\Omega)} \le M\|v\|_{L^2(\Omega)}\|z_{u,v}\|_{L^2(\Omega)} \le M\|v\|_{L^2(\Omega)}\|z_{u,v}\|_{H^1(\Omega)}.
\]
Hence, we have $\|z_{u,v}\|_{H^1(\Omega)} \le M[\Lambda_{\bar u} - C^2_{4,\Omega}\varepsilon_{\bar u}]^{-1}\|v\|_{L^2(\Omega)}$ for every $u \in B_{\varepsilon_{\bar u}}(\bar u)$. Using this fact and \cite{Nittka11}, the estimate in $H^1(\Omega) \cap C^{0,\mu}(\bar\Omega)$ uniformly in the ball $B_{\varepsilon_{\bar u}}(\bar u)$ follows. Therefore, we get that $\|G'(u)\| \le L_G$ for all $u \in B_{\varepsilon_{\bar u}}(\bar u)$ and some constant $L_G$.

The proof of the Lipschitz property of $\Phi$ in $B_{\varepsilon_{\bar u}}(\bar u)$ follows the same steps as for $G$. First, we prove that $\|\varphi_u\|_{H^1(\Omega) \cap C^{0,\mu}(\bar\Omega)} \le K$ for all $u \in B_{\varepsilon_{\bar u}}(\bar u)$ and some constant $K$. To this end, we use that $\|y_u\|_{C(\bar\Omega)} \le M$. Hence, with Assumption \ref{A3.2} we infer that $\Big|\frac{\partial L}{\partial y}(x,y_u(x))\Big| \leq L_M(x)$, where $L_M \in L^p(\Omega)$. Therefore, the inequality $\|\varphi_u\|_{H^1(\Omega) \cap C^{0,\mu}(\bar\Omega)} \le K$ follows from the equation \eqref{E3.3}. Then, we apply the mean value theorem and estimate $\|\eta_{u,v}\|_{H^1(\Omega) \cap C^{0,\mu}(\bar\Omega)} = \|\Phi'(u)v\|_{H^1(\Omega) \cap C^{0,\mu}(\bar\Omega)}$ for all $\|v\|_{L^2(\Omega)} \le 1$ and $u \in B_{\varepsilon_{\bar u}}(\bar u)$. This is obtained from equation \eqref{E3.4} using the $C(\bar\Omega)$ estimates for $y_u$ and $\varphi_u$ as well as the $H^1(\Omega)$ estimates for $z_{u,v}$. Finally, we get the Lipschitz property for $S$:
\begin{align}
&\|S(u_2) - S(u_1)\|_{H^1(\Omega) \cap C^{0,\mu}(\bar\Omega)}\notag\\
& \le \frac{1}{\tichonov}\Big[\|G(u_2) - G(u_1)\|_{H^1(\Omega) \cap C^{0,\mu}(\bar\Omega)}\|\Phi(u_2)\|_{H^1(\Omega) \cap C^{0,\mu}(\bar\Omega)}\notag\\
& + \|\Phi(u_2) - \Phi(u_1)\|_{H^1(\Omega) \cap C^{0,\mu}(\bar\Omega)}\|G(u_1)\|_{H^1(\Omega) \cap C^{0,\mu}(\bar\Omega)}\Big]\notag\\
& \le L_S\|u_2 - u_1\|_{L^2(\Omega)}\quad \forall u_1, u_2 \in B_{\varepsilon_{\bar u}}(\bar u).
\label{E4.3}
\end{align}
This completes the proof of the $\partial\Psi$-semismoothness of $\Psi$ in $\A$.\qed
\end{proof}
\begin{corollary}
The function $F:\A \longrightarrow L^2(\Omega)$ is $\partial F$-semismooth in $\A$, where
\[
\partial F(u) = \{M = I - N : N \in \partial\Psi(u)\}.
\]
\label{C4.2}
\end{corollary}
This is a straightforward consequence of Theorem \ref{T4.3}.

To implement Algorithm 1, we select the operators $M_u:L^2(\Omega) \longrightarrow L^2(\Omega)$ for every $u \in \A$ as follows. First, we define the function $\lambda:\mathbb{R} \longrightarrow \mathbb{R}$ by
\[
\lambda(t) = \Big\{\begin{array}{cl} 1&\text{if } t \in (\umin,\umax),\\0&\text{otherwise}.\end{array}
\]
It is obvious that $\lambda(t) \in \partial\psi(t)$ for every $t \in \mathbb{R}$. We define $M_u:L^2(\Omega) \longrightarrow L^2(\Omega)$ by $M_uv = v - h_u\cdot S'(u)v$, where $h_u(x) = \lambda(S(u)(x)) = \lambda\big(\frac{1}{\tichonov}y_u(x)\varphi_u(x)\big)$. It is immediate that $M_u \in \partial F(u)$. For this selection we have the following result.

\begin{theorem}
Let $(\bar u,\bar y,\bar\varphi) \in \uad \times [H^1(\Omega) \cap C^{0,\mu}(\bar\Omega)]^2$ satisfy the first order optimality conditions \eqref{E3.5}--\eqref{E3.7}, the strict complementarity condition $|\Sigma_{\bar u}| = 0$, and the second order sufficient optimality condition $J''(\bar u)v^2 > 0$ for every $v \in C_{\bar u}\setminus\{0\}$. Then, there exist $\delta > 0$ and $C > 0$ such that for every $u \in B_\delta(\bar u) \subset \A$ the linear operator $M_u:L^2(\Omega) \longrightarrow L^2(\Omega)$ is an isomorphism and $\|M_u^{-1}\| \le C$ holds.
\label{T4.4}
\end{theorem}

\begin{proof}
Given $u \in \A$ we define the active and inactive sets for $u$ as follows
\[
\mathbb{A}_u = \{x \in \Omega : \frac{1}{\tichonov}y_u(x)\varphi_u(x) \not\in (\umin,\umax)\},\]
\[ \mathbb{I}_u = \{x \in \Omega : \frac{1}{\tichonov}y_u(x)\varphi_u(x) \in (\umin,\umax)\}.
\]
We denote by $\chi_{\mathbb{A}_u}$ and $\chi_{\mathbb{I}_u}$ the characteristic functions of $\mathbb{A}_u$ and $\mathbb{I}_u$, respectively.
According to the definition of $M_u$ we have $M_uv = v- \frac{1}{\tichonov} [z_{u,v}\varphi_u + y_u \eta_{u,v}]\ci$. It is obvious that $M_u$ is a continuous operator. Let us prove that for every $w \in L^2(\Omega)$ there exists a unique $v \in L^2(\Omega)$ such that $M_uv = w$. This equation can be written in the form
\[
\left\{\begin{array}{ll} v(x) = w(x) &\text{if } x \in \mathbb{A}_u,\\v(x) - \frac{1}{\tichonov} [z_{u,v}(x)\varphi_u(x) + y_u(x)\eta_{u,v}(x)] = w(x)&\text{if } x \in \mathbb{I}_u.\end{array}\right.
\]
Taking into account that $v$ coincides with $w$ in $\mathbb{A}_u$ and, hence, $v = \chi_{{}_{\mathbb{I}_u}}v + \chi_{{}_{\mathbb{A}_u}}w$, the second equation can be written
\begin{equation}
\ci v - \frac{1}{\tichonov} [z_{u,\ci v}\varphi_u + y_u\eta_{u,\ci v}] = w + \frac{1}{\tichonov} [z_{u,\ca w}\varphi_u + y_u\eta_{u,\ca w}].
\label{E4.4}
\end{equation}
In order to prove the existence and uniqueness of a solution of \eqref{E4.4} we introduce the quadratic functional $\mathbb{J}:L^2(\mathbb{I}_u) \longrightarrow \mathbb{R}$ defined by
\begin{align*}
\mathbb{J}(v) &= \frac{1}{2\tichonov}J''(u)(\ci v)^2 - \int_{\mathbb{I}_u}\big(w + \frac{1}{\tichonov} [z_{u,\ca w}\varphi_u + y_u\eta_{u,\ca w}])v\dx\\
&=\frac{1}{2}\int_{\mathbb{I}_u}\big[v^2 -\frac{1}{\tichonov}(\varphi_uz_{u\ci v} + y_u\eta_{u,\ci v})v\big]\dx\\
& - \int_{\mathbb{I}_u}\big(w + \frac{1}{\tichonov} [z_{u,\ca w}\varphi_u + y_u\eta_{u,\ca w}])v\dx,
\end{align*}
where we have used the expression of $J''(u)$ given in Corollary \ref{C3.1}. We observe that $\mathbb{J}'(v) = 0$ if and only if $v$ satisfies \eqref{E4.4}. Therefore, if we prove that $\mathbb{J}$ has a unique stationary point, then the existence and uniqueness of a solution of \eqref{E4.4} follows. From Theorem \ref{T3.5} we get that \eqref{E3.12} holds for some $\tau > 0$ and $\kappa > 0$. Since $J''$ is a continuous functional in $\A$, we deduce the existence of $\delta_0 > 0$ such that $|[J''(u) - J''(\bar u)]v^2| \le \frac{\kappa}{2}\|v\|^2_{L^2(\Omega)}$ for all $v \in L^2(\Omega)$ if $\|u - \bar u\|_{L^2(\Omega)} \le \delta_0$. This inequality and \eqref{E3.12} imply that
\begin{equation}
J''(u)v^2 \ge \frac{\kappa}{2}\|v\|^2_{L^2(\Omega)}\quad \forall v \in T^\tau_{\bar u}\text{ and } \forall u \in B_{\delta_0}(\bar u).
\label{E4.5}
\end {equation}
Now, we prove that $\ci v \in T^\tau_{\bar u}$ for all $v \in L^2(\mathbb{I}_u)$ if $u$ is sufficiently close to $\bar u$. As a consequence we infer that the quadratic form $\mathbb{J}$ is strictly convex and coercive, which proves the existence of a unique stationary point, the unique minimizer. To prove that $\ci v \in T^\tau_{\bar u}$ we select $\delta = \min\{\delta_0,\varepsilon_{\bar u},\frac{\tau}{\tichonov L_S}\}$, where $L_S$ was given in \eqref{E4.3}. Hence, we have that $\|S(u) - S(\bar u)\|_{C(\bar\Omega)} \le L_S\delta \le \frac{\tau}{\tichonov}$ for every $u \in B_\delta(\bar u)$. If $\tichonov\bar u(x) - \bar y(x)\bar\varphi(x) > \tau$, then
\eqref{E3.7} implies that $\bar u(x) = \umin$ and, hence, $\frac{1}{\tichonov}\bar y(x)\bar\varphi(x) < \umin - \frac{\tau}{\tichonov}$. This yields $\frac{1}{\tichonov}y_u(x)\varphi_u(x) < \umin$, therefore we have $x \in \mathbb{A}_u$ and $(\ci v)(x) = 0$. Analogously we proceed if $\tichonov\bar u(x) - \bar y(x)\bar\varphi(x) < -\tau$.

It remains to deduce the existence of a constant $C$ such that $\|M_u^{-1}\| \le C$ for every $u \in B_\delta(\bar u)$. From \eqref{E4.5}, Corollary \ref{C3.1}, and \eqref{E4.4} we infer
\begin{align*}
\frac{\kappa}{2}\|\ci v\|^2_{L^2(\Omega)} &\le J''(u)(\ci v)^2 = \int_\Omega(\tichonov v -[\varphi_uz_{u,\ci v} + y_u\eta_{u,\ci v}])\ci v\dx\\
&= \int_\Omega(\tichonov w + [\varphi_uz_{u,\ca w} + y_u\eta_{u,\ca w}])\ci v\dx.
\end{align*}
Since $\ca w = \ca v$ we have
\[
\tichonov\|\ca v\|^2_{L^2(\Omega)} = \tichonov\int_\Omega w\ca v\dx.
\]
From the last to relations, the uniform estimates for $y_u$ and $\varphi_u$ in the ball $B_{\varepsilon_u}(\bar u)$ established in the proof of Theorem \ref{T4.3}, and equations \eqref{E2.9} and \eqref{E3.4} we get the existence of a constant $C'$ such that
\begin{align*}
\min\big\{\tichonov,\frac{\kappa}{2}\big\}\|v\|^2_{L^2(\Omega)}&\le \tichonov\int_\Omega w v\dx + \int_\Omega[\varphi_uz_{u,\ca w} + y_u\eta_{u,\ca w}]\ci v\dx\\
&\le C'\|w\|_{L^2(\Omega)}\|v\|_{L^2(\Omega)},
\end{align*}
which proves that $\|M_u^{-1}\| \le \frac{C}{\min\{\tichonov,\frac{\kappa}{2}\}}$.
\qed\end{proof}

Algorithm  \ref{Alg2} implements the semismooth Newton method to solve \Pb.
 The following corollary establishes its convergence.

\begin{algorithm2e}[H]
\caption{Semismooth Newton method to solve \Pb.}\label{Alg2}
\DontPrintSemicolon
Initialize. Choose $u_0\in L^2(\Omega)$. Set $j=0$.\;
\For{$j\geq 0$}{
Compute $y_j = G(u_j)$ solving the nonlinear equation\label{line3alg2}
\[A y + a(x,y)+ u_jy =0\text{ in }\Omega,\ \partial_{\tichonov_{A}} y=g\text{ in }\Gamma\]
\;
Compute $\varphi_j = \Phi(u_j)$ solving the linear equation
\[ A^*\varphi + \displaystyle \frac{\partial a}{\partial y}(x,y_j)\varphi + u_j\varphi = \frac{\partial L}{\partial y}(x,y_j) \   \mbox{in } \Omega,\  \partial_{\conormal_{A^*}} \varphi = 0\  \mbox{on }\Gamma.
\]
\;
Set $\mathbb{A}_j =\mathbb{A}^{+}_j\cup\mathbb{A}^{-}_j$ and  $\mathbb{I}_j=\Omega\setminus\mathbb{A}_j$, where
\[\mathbb{A}^{+}_j = \{x\in\Omega: y_j(x)\varphi_j(x)\geq \tichonov\umax\},\]
\[\mathbb{A}^{-}_j = \{x\in\Omega: y_j(x)\varphi_j(x)\leq \tichonov\umin\}\]
\;
Set
\[w_j(x) = -F(u_j)(x) = \left\{\begin{array}{ll}
                              -u_j(x) +\umax & \text{ if }x \in \mathbb{A}^{+}_j \\
                              -u_j(x) + \frac{1}{\tichonov}\varphi_j(x)y_j(x) & \text{ if }x \in \mathbb{I}_j \\
                              -u_j(x) +\umin & \text{ if }x \in \mathbb{A}^{-}_j
                            \end{array}\right.
\]
\;
Compute $z_j = z_{u_j,\chi_{{}_{\mathbb{A}_j}}w_{j}}$ and $\eta_j = \eta_{u_j,\chi_{{}_{\mathbb{A}_j}}w_{j}}$ solving the linear equations\label{line8}
\begin{align*}
A z_j +\frac{\partial a}{\partial y}(x,y_j)z_j + u_j z_j =& -y_j \chi_{{}_{\mathbb{A}_j}}w_{j}\text{ in }\Omega,\ \partial_{\tichonov_{A}} z_j = 0\text{ on }\Gamma\\
A^*\eta_j +\frac{\partial a}{\partial y}(x,y_j)\eta_j +u_j \eta_j = & \left(
\frac{\partial^2 L}{\partial y^2}(x,y_j)
-\varphi_j \frac{\partial^2 f}{\partial y^2}(x,y_j)\right)z_j\\
& -\varphi_j \chi_{{}_{\mathbb{A}_j}}w_{j}\text{ in }\Omega,\, \partial_{\conormal_{A^*}} \eta_j = 0\text{ on }\Gamma
\end{align*}
\;
Solve the quadratic problem\label{line9}
\[
\mathrm{(}Q_j\mathrm{)}\qquad \min_{v\in L^2(\mathbb{I}_j)} \mathbb{J}_j(v) := \frac{1}{2\tichonov} J''(u_j)(\chi_{{}_{\mathbb{I}_j}} v)^2 -\int_{\mathbb{I}_j}( w_j +\frac{1}{\tichonov}[z_j\varphi_j+y_j\eta_j]) v\mathrm{d}x
\]
Name $v_{\mathbb{I}_j}$ its solution.\;
Set $v_{j} =  \chi_{{}_{\mathbb{A}_j}}w_{j} + \chi_{{}_{\mathbb{I}_j}} v_{\mathbb{I}_j}$\;
Set $u_{j+1}=u_j+v_j$ and $j=j+1$.\;
}
\end{algorithm2e}

\begin{corollary}
Let $(\bar u,\bar y,\bar\varphi) \in \uad \times [H^1(\Omega) \cap C^{0,\mu}(\bar\Omega)]^2$ satisfy the first order optimality conditions \eqref{E3.5}--\eqref{E3.7}, the strict complementarity condition $|\Sigma_{\bar u}| = 0$, and the second order sufficient optimality condition $J''(\bar u)v^2 > 0$ for every $v \in C_{\bar u}\setminus\{0\}$. Then, there exists $\delta > 0$ such that for all $u_0 \in B_\delta(\bar u)$ the sequence $\{u_j\}$ generated by Algorithm \ref{Alg2} is contained in the ball $B_\delta(\bar u)$ and converges superlinearly to $\bar u$.
\label{C4.3}
\end{corollary}

\begin{proof}
Since any local solution of \Pb satisfying the second order sufficient condition is locally the unique stationary point of \Pb, see \cite{Casas-Troltzsch12}, this result is a straightforward consequence of Theorem \ref{T4.1}, Corollary \ref{C4.2}, and Theorem \ref{T4.4}.
\qed\end{proof}

\begin{remark}
To solve the quadratic problem $\mathrm{(}Q_j\mathrm{)}$ that appears in line \ref{line9} of Algorithm \ref{Alg2} we can use, e.g., the conjugate gradient method. Notice that we can write $\mathbb{J}_j(v) = \frac{1}{2}(v,A_j v)_{L^2(\mathbb{I}_j)} - (b_j,v)_{L^2(\mathbb{I}_j)}$, where $b_j = \chi_{{}_{\mathbb{I}_j}} ( w_{j} +\frac{1}{\tichonov}[z_j\varphi_j+y_j\eta_{j}])$ and we can compute $A_j v$ using Algorithm \ref{Alg3}. Therefore $\mathrm{(}Q_j\mathrm{)}$ can be solved without need of the explicit computation of the Hessian $J''(u_j)$.
\end{remark}
\LinesNumbered
\begin{algorithm2e}[H]
\caption{Computation of the product Hessian-vector}\label{Alg3}
\DontPrintSemicolon
Solve
\begin{align*}
A z +\frac{\partial a}{\partial y}(x,y_j)z + u_j z =& -y_j \chi_{{}_{\mathbb{I}_j}}v\text{ in }\Omega,\ \partial_{\tichonov_{A}} z = 0\text{ on }\Gamma
\end{align*}
\;
Solve
\begin{align*}
A^*\eta +\frac{\partial a}{\partial y}(x,y_j)\eta +u_j \eta = & \left(
\frac{\partial^2 L}{\partial y^2}(x,y_j)
-\varphi_j \frac{\partial^2 f}{\partial y^2}(x,y_j)\right)z\\
& -\varphi_j \chi_{{}_{\mathbb{I}_j}}v\text{ in }\Omega,\, \partial_{\conormal_{A^*}} \eta = 0\text{ on }\Gamma
\end{align*}
\;
Set $A_jv = \chi_{{}_{\mathbb{I}_j}}(v+\frac{1}{\tichonov}[z\varphi_j+\eta y_j])$\;
\end{algorithm2e}

\section{A numerical example}
We have used Algorithm \ref{Alg2} to solve the problem with the following data: $\Omega=(0,1)^2\subset\mathbb{R}^2$, $A=-\Delta$, $g(x)=0$,
\[f(x,y) = y^3|y|+2y-100\sin(2\pi x_1)\sin(\pi x_2),\]
$\tichonov = 0.05$, $\umin = -1$, $\umax = 1$, $L(x,y) = 0.5(y-y_d(x))$, where
\[y_d(x) = -64 x_1(1-x_1) x_2(1-x_2).\]
We solve a finite element discretization of \Pb. Continuous piecewise linear functions are used for the state, the adjoint state, and the control. The Tichonov regularization term is discretized using the lump mass matrix. In this way, the optimization algorithm for the discrete problem is exactly the discrete version of Algorithm \ref{Alg2}.

The chosen initial point is $u_0=0$. The algorithm stops when \[\delta_j = \frac{\| v_j\|_{L^2(\Omega)}}{\max\{1,\| u_{j+1}\|_{L^2(\Omega)}\}} <  5\times 10^{-14}\] or when $J(u_j)$ and $J(u_{j+1})$ are equal up to machine precision. At each iteration, the solution of the quadratic subproblem $(Q_j)$ is obtained using the conjugate gradient method implemented in Matlab built-in command \texttt{pcg} and the nonlinear equation in line \ref{line3alg2} is solved using Newton's method. The tolerance  $5\times 10^{-14}$ is used for both subproblems.

We show the convergence results for different mesh sizes in tables \ref{T1}--\ref{T3}. Not only the predicted superlinear convergence can be observed in all of them, but also the mesh-independence of the convergence history, which is explained thanks to the convergence result for the algorithm in the infinite-dimensional setting.

$\sharp$Newton is the number of Newton iterations to solve the nonlinear PDE in line \ref{line3alg2} and $\sharp$CG is the number of iterations of the conjugate gradient method used to solve $(Q_j)$ in Algorithm \ref{Alg2}. In terms of computational cost, the hard work is done to solve the nonlinear PDE: at each of the $\sharp$Newton iterations we have to factorize a (sparse) matrix. The last factorization obtained at this step can be used to solve the rest of the linear PDEs that appear in Algorithm \ref{Alg2} and also the ones involved in the conjugate gradient method, so a good measure of the computational cost is given by the total amount of $\sharp$Newton iterations.

\begin{table}[h]
  \centering
  \begin{tabular}{ccccc}
    $j$ &       $J(u_j)$            &   $\delta_j$   &  $\sharp$Newton & $\sharp$CG  \\ \hline
      0&  3.9142466314434916e+00 &   8.8e-01  &   8 &  13   \\
  1&  3.8210815158943565e+00 &   8.8e-03  &   5 &  12  \\
  2&  3.8210805974920747e+00 &   3.5e-09  &   3 &  13  \\
  3&  3.8210805974920712e+00 &   5.9e-15  &   2 &  14  \\
  4 & 3.8210805974920659e+00 &            &  1    &
  \end{tabular}
  \caption{Convergence history of the problem in the example for $h=2^{-7}$. }\label{T1}
\end{table}

\begin{table}[h]
  \centering
  \begin{tabular}{ccccc}
    $j$ &       $J(u_j)$            &   $\delta_j$   &  $\sharp$Newton & $\sharp$CG  \\ \hline
  0&  3.9149225191422081e+00 &   8.8e-01  &   8 &  12    \\
  1&  3.8217486072447922e+00 &   9.4e-03  &   5 &  12    \\
  2&  3.8217477897599132e+00 &   4.4e-07  &   3 &  12    \\
  3&  3.8217477897599528e+00 &   3.7e-14  &   2 &  12    \\
  4 & 3.8217477897599559e+00 &            &  1    &
  \end{tabular}
  \caption{Convergence history of the problem in the example for $h=2^{-8}$.}\label{T2}
\end{table}

\begin{table}[h]
  \centering
  \begin{tabular}{ccccc}
    $j$ &       $J(u_j)$            &   $\delta_j$   &  $\sharp$Newton & $\sharp$CG  \\ \hline
  0&  3.9150915018042891e+00 &   8.8e-01  &   8 &  12 \\
  1&  3.8219154943064222e+00 &   9.6e-03  &   5 &  12 \\
  2&  3.8219145854336314e+00 &   4.1e-07  &   3 &  12 \\
  3&  3.8219145854337437e+00 &   3.6e-14  &   2 &  13 \\
  4 & 3.8219145854337260e+00 &            &  1  &
  \end{tabular}
  \caption{Convergence history of the problem in the example for $h=2^{-9}$.}\label{T3}
\end{table}

\begin{figure}[h]
 \centering
 \includegraphics[width=.5\textwidth]{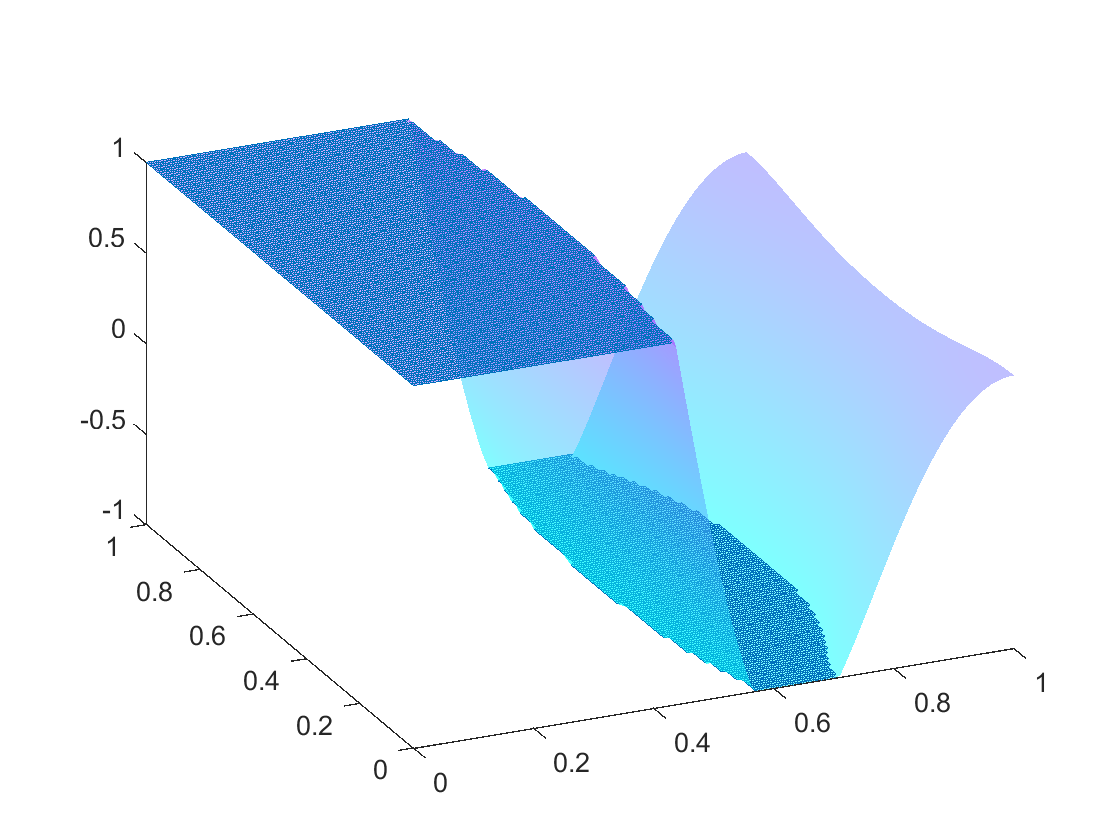}
   \caption{Optimal control for the example.}\label{F1}
\end{figure}

A picture of the optimal control can be seen in Figure \ref{F1}. For the finest mesh, we find that $|\{x\in\Omega:\ \bar u(x)=\umax\}| = 0.459$, $|\{x\in\Omega:\ \bar u(x)=\umin\}| = 0.233$, $|\{x\in \Omega:\ \umin < \bar u(x) < \umax\}| = 0.308$ and $|\Sigma_{\bar u}|=0$.

\end{document}